# Effects of End Restraints on Thermo-Mechanical Response of Energy Piles

Arash Saeidi Rashk Olia[1], Dunja Perić[2]

[1] Department of Civil Engineering, Kansas State University, 1701C Platt St., Manhattan, KS 66506-5000; e-mail: saeidi@ksu.edu

[2] Department of Civil Engineering, Kansas State University, 1701C Platt St., Manhattan, KS, 66506-5000; e-mail: peric@ksu.edu

## ABSTRACT

Currently, soil structure interaction in energy piles has not been understood thoroughly. One of the important underlying features is the effect of tip and head restraints on displacement, strain and stress in energy piles. This study has investigated thermo-mechanical response of energy piles subjected to different end restraints by using recently found analytical solutions, thus providing a fundamental, rational, mechanics-based understanding. End restraints are found to have a substantial effect on thermo-mechanical response of energy piles, especially on thermal axial displacement and axial stress in the pile. Head restraint imposed by interaction of an energy pile with the superstructure led to a decrease in the magnitude of head displacement and increase in axial stress, while decreasing the axial strain. The impact of head restraint was more pronounced in end bearing than in fully floating energy piles.

## INTRODUCTION

Energy piles are dual purpose deep foundations. They enable the exchange of thermal energy between the superstructure and shallow subsurface while transferring the superstructure loads to the ground. Thus, they use naturally renewable subsurface geothermal energy as a supplemental energy source for space heating and cooling (Brandl 2006, Ghasemi-fare and Basu 2015).

While the energy exchange between energy piles and the surrounding soil induces the temperature difference to the pile element, free thermal deformations in energy piles is restrained by the surrounding soil, superstructure and soil layers beneath the pile tip. Therefore, once the pile is heated, and has the tendency to expand these restraints induce compressive stresses to the pile element and vice versa, as the energy pile is cooled and tends to contract tensile stress is generated in the pile. Consequently, design of energy piles requires a comprehensive understanding of the both thermal and mechanical responses of the pile, as well as the resultant combined thermo-mechanical behavior of these structural elements.

Schematic diagrams of thermo-mechanical response of a single energy pile subjected to thermal and mechanical loads were presented by Bourne-Webb et al. (2009), Amatya et al. (2012), and Bourne-Webb et al. (2013). They presented thermal axial strain, load and shaft shear stress for different end restraints based on the results of full scale in situ tests. More recently Rotta Loria and Laloui (2018) provided axial stress, displacement and soil pile interface shear stress in an energy pile subjected to thermo-mechanical loading. While they did not present any governing equations, other than those enforcing global equilibrium, they stated that their load transfer diagrams were based on thermo-elasticity theory. Saeidi Rashk Olia and Perić conducted numerical modeling using finite element analysis to evaluate thermo-mechanical response of energy piles for semi-floating and end bearing piles. The results of numerical models were successfully validated against analytical solutions used herein (Saeidi Rashk Olia and Perić 2021 a & b).

The current study investigates effects of head and end restraints on combined thermo-mechanical axial displacement, strain and stress in a single energy pile that is embedded in a homogenous soil, based on the analytical solutions.

**ANALYTICAL MODEL**

Analytical solutions for thermo-mechanical soil structure interaction of a single energy pile were presented by Cossel (2019), Iodice et al. (2020) Perić et al. (2020). The solutions were derived based on the assumption that the energy pile obeys the thermo-elastic constitutive law while the soil pile interface remains elastic. These assumptions were confirmed by Knellwolf et al. (2011) and Perić et. al. (2017). The physical model depicted in Figure 1 accounts for restraining effects of the surrounding soil. Specifically, a continuous linear elastic shear spring having stiffness of $k_s$ represents restraint imposed by surrounding soil to free deformation in energy piles, while linear elastic normal springs of stiffness $k_h$ and $k_b$ represent head and tip restraints imposed by the superstructure and soil layers respectively.

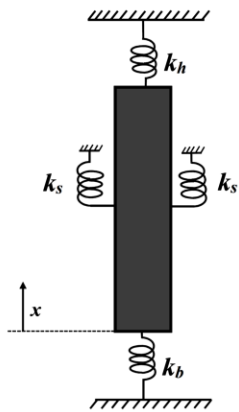

**Figure 1. Schematic of the physical model of energy pile**

Cossel (2019) and Perić, et al. (2020) provided the analytical solutions for axial displacement, strain and stress in an ideally end bearing energy pile (EB) subjected to thermal and mechanical loading. The solutions are given by Eqs. 1, 2 and 3 respectively. These solutions were derived for mechanical and thermal loading by imposing equilibrium to the infinitesimal pile segment in vertical direction. These derivations were based on the assumption that the energy pile follows thermo-elastic constitutive law, as well as that the soil pile interface remains in elastic range. Furthermore, the shear stress at the interface is a linear function of the vertical displacement. Temperature change was also assumed to be constant throughout the entire pile.

$$u(x) = \frac{\alpha \Delta T \sinh(\psi x)}{\psi \cosh(\psi L) + \frac{k_h}{E}\sinh(\psi L)} + \frac{F \sinh(\psi x)}{A E \psi \cosh(\psi L)} \quad (1)$$

$$\varepsilon(x) = \frac{\alpha \Delta T \cosh(\psi x)}{\cosh(\psi L) + \frac{k_h}{E\psi}\sinh(\psi L)} + \frac{F \cosh(\psi x)}{A E \cosh(\psi L)} \quad (2)$$

$$\sigma(x) = E\alpha \Delta T \left[\frac{\cosh(\psi x)}{\cosh(\psi L) + \frac{k_h}{E\psi}\sinh(\psi L)} - 1\right] + \frac{F \cosh(\psi x)}{A \cosh(\psi L)} \quad (3)$$

where $x$ is the vertical coordinate with positive direction pointing upward as shown in Figure 1, and $u$, $\varepsilon$ and $\sigma$ are axial displacement, strain and stress in the energy pile respectively. Positive displacement points upward, while tensile strain and stress are positive. The elastic modulus and the coefficient of thermal expansion of the pile are denoted by $E$ and $\alpha$ respectively, while $\Delta T$ is a temperature difference of the pile relative to the surrounding soil. Heating is represented by positive $\Delta T$ while cooling implies negative $\Delta T$. Mechanical load in the form of an axial force $F$ is positive in tension and negative in compression. The parameter $\psi$ contains the pile geometry and relative stiffness of the soil with respect to the pile. It is a given by

$$\psi^2 = \left(\frac{p}{A}\right)\left(\frac{k_s}{E}\right) \quad (4)$$

where $p$ and $A$ are the perimeter and the cross-sectional area of the energy pile respectively. The length of the pile is denoted by $L$.

In case of an ideally fully floating (FF) energy piles whereby no tip restraint is present, the solutions for displacement, strain and stress in a semi floating pile subjected to thermo-mechanical loading presented by Cossel (2019) are modified by setting $k_b \to 0$, thus giving

$$u(x) = \frac{\alpha \Delta T \sinh[\psi(x-x_0)]}{\psi \cosh[\psi(L-x_0)] + \frac{k_h}{E}\sinh[\psi(L-x_0)]} + \frac{F \cosh(\psi x)}{A E \psi \sinh(\psi L)} \quad (5)$$

$$\varepsilon(x) = \frac{\alpha \Delta T \cosh[\psi(x-x_0)]}{\cosh\psi[(L-x_0)] + \frac{k_h}{E\psi}\sinh[\psi(L-x_0)]} + \frac{F \sinh(\psi x)}{A E \sinh(\psi L)} \quad (6)$$

$$\sigma(x) = E\alpha \Delta T \left[\frac{\cosh[\psi(x-x_0)]}{\cosh\psi[(L-x_0)] + \frac{k_h}{E\psi}\sinh[\psi(L-x_0)]} - 1\right] + \frac{F \sinh(\psi x)}{A \sinh(\psi L)} \quad (7)$$

where $x_0$ is the location of the zero thermal displacement or thermal null point. In ideally end bearing energy piles thermal null point is located at the pile tip regardless of the amount of a head restraint. Location of a thermal null point in a fully floating energy pile is obtained from Cossel (2019) by simply setting $k_b \to 0$. It results in

$$x_0 = \frac{1}{\psi} tanh^{-1} \left[ \frac{cosh(\psi L) - 1 + \frac{k_h}{E\psi} sinh(\psi L)}{sinh(\psi L) + \frac{k_h}{E\psi} cosh(\psi L)} \right] \tag{8}$$

Table. 1 summarizes parameters of analytical solutions presented by Cossel (2019) and Perić, et al. (2020)

Table 1. Parameters of the analytical model

| Symbol | Parameter | Unit |
|---|---|---|
| $\varphi$ | Geometry and stiffness ratio parameter | 1/m |
| $k_s$ | Shear spring stiffness | GPa/m |
| $k_h$ | Head spring stiffness | GPa/m |
| $L$ | Energy pile length | m |
| $\alpha$ | Coefficient of thermal expansion | 1/°C |
| $E$ | Pile element elastic modulus | GPa |
| $\Delta T$ | Temperature difference of energy pile with surrounding soil | °C |
| $F$ | Mechanical load applied to the pile head | kN |

**MATERIAL PROPERTIES**

For the purpose of thermo-mechanical analyses conducted herein, the energy pile is embedded in a single soil layer, the bottom of which is flush with the pile tip. The layer is underlain by the bedrock. To this end, one of the four soil layers surrounding the energy pile constructed and tested at Swiss Federal Institute of Technology Lausanne, Switzerland (Laloui et al. 2006) was selected to surround the energy pile herein. Based on laboratory classification tests, the results of which were provided by Laloui et al. (1999), Perić et al. (2017) classified all four layers as low plasticity clays (CL). In the present study, the shear spring stiffness ($k_s$) was chosen to be equal to 0.0167 GPa/m as reported by Knellwolf et al. (2011) for soil layer A1 in Swiss Federal Institute of Technology Lausanne. It is noted that only two cases of tip restraints are considered herein. They correspond to an ideally end bearing and ideally fully floating energy piles, thus resulting in infinitely large and zero $k_b$ values respectively. While realistically the stiffness always has a finite value that is never zero these two cases represent limits to any actual situation. Although the length of an ideally fully floating pile would be infinite, a finite length can be combined with an arbitrary small value of $k_b$ that corresponds to a desired arbitrary small ratio of the mechanically induced stresses at the pile tip and pile head. The resulting value of $k_b$ then corresponds to that of a nearly fully floating pile, which is simply referred to as a fully floating pile herein. Similarly, the ideally end bearing pile is referred to as simply end bearing pile.

Laloui et al. (2006) reported the actual pile to be 26 m long, its diameter is equal to 1 m, the pile element's modulus of elasticity (*E)* was also reported to be 29.2 GPa, and its coefficient of the thermal expansion (*α*) is equal to $1 \times 10^{-5}$ /°C.

## LOAD SCENARIOS

Two thermo-mechanical load scenarios that involve compressive mechanical load are considered. They are described as follows:

- Load scenario (*i*): $F < 0$ (compression) and $\Delta T < 0$ (cooling)
- Load scenario (*ii*): $F < 0$ (compression) and $\Delta T > 0$ (heating)

It is noted that in the case of scenario (*i*) the pile stress induced by thermal loading has a different sign from the one induced by mechanical loading. On the contrary, in the case of loading scenario (*ii*) thermal load induces the stress of the same sign as the one induced by mechanical load. Thus, in the latter case the pile stress has the same sign along the entire length of the pile regardless of a magnitude of thermal and mechanical loads.

Selected magnitude of compressive load is 1000 kN, thus resembling the 4-storey building at Swiss Federal Institute of Technology (Laloui et al 2006). Furthermore, thermal load of $\Delta T = \pm 10°C$ was applied in both cases, single floating and single end bearing energy piles.

Perić et al. (2020) defined equivalent thermal load $\Delta T_{eq}$, thus enabling a direct conversion between thermal and mechanical loads. It is given by

$$\Delta T_{eq} = \pm \frac{F}{AE\alpha} \tag{9}$$

where $A$ is the cross sectional area of the pile, and plus sign applies to load scenario (*i*) while minus sign applies to the scenario (*ii*). For a mechanical load of magnitude $|F|$, the corresponding magnitude of the equivalent thermal load is $|\Delta T_{eq}|$. Eq. 9 is based on the fact mechanical load produces displacement and strain of equal magnitude as the one produced by the equivalent thermal load, throughout the entire length of an end bearing pile. Thus, for the mechanical load $F=-1000$ kN Eq. 9, results in the magnitude of equivalent thermal load being equal to 4.36 °C. Consequently, for the selected thermal load of $\pm 10$ °C the magnitude of actual thermal load is 2.29 times larger than that of the equivalent thermal load whereby the actual mechanical load is -1000 kN.

## RESULTS

Thermal, mechanical and combined responses to load scenario (*i*) are depicted in Figures 2 and 3 for fully floating and end bearing piles respectively. Both of these figures indicate that magnitudes of thermally induced displacement and strain are larger than those of mechanically induced displacement and strain. So, the overall displacement and strain responses are dominated by the thermal response. In the case of the stress response, it is seen from Figures 2 and 3 that thermal stress becomes dominant at certain depth, thus resulting in development of tensile stress and tension zone in the pile. The tensile stresses are larger and tension zone is longer in the end bearing pile than in the fully floating pile. Clearly, if magnitude of thermal load were sufficiently small no tensile stress would have developed.

Otherwise, magnitude of the combined displacement at the pile head is larger in the end bearing pile than in the fully floating pile and vice versa at the pile tip. This is the consequence of the fact that the thermal null point is located at the mid length of the pile in fully floating pile and at the pile tip in end bearing pile. Furthermore, the magnitude of a combined strain is larger in the fully floating pile than in the end bearing pile throughout the length of the pile. As expected, the magnitude of stress and strain at the pile tip of a fully floating pile are nearly equal to zero. Also, in the absence of a head restraint the magnitude of thermally induced stress at the pile head is equal to zero.

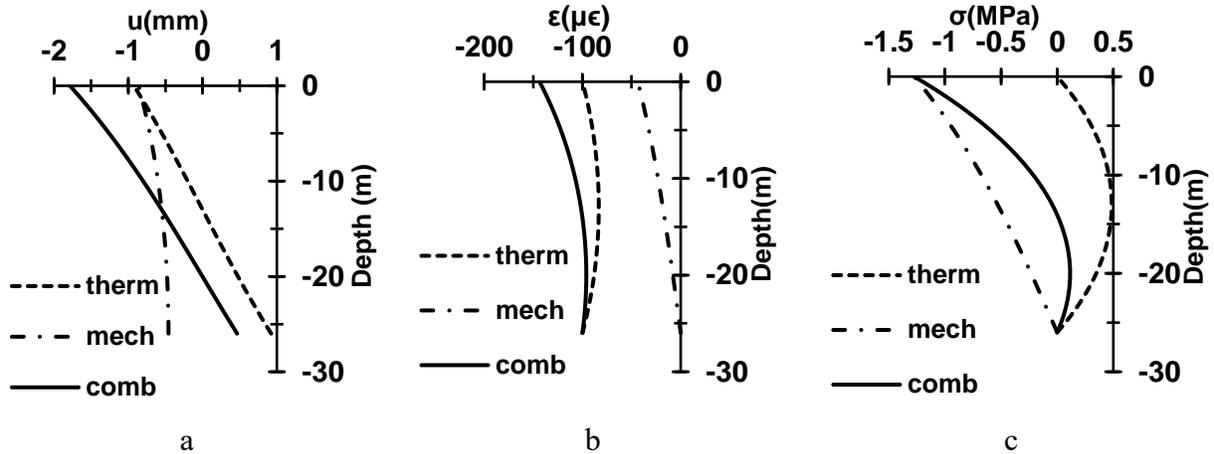

**Figure 2. a) Displacement, (b) strain, and (c) stress in a fully floating energy pile subjected to thermo-mechanical load scenario (*i*), $F = -1\text{M N}$ and $\Delta T = -10$ °C, $k_h = 0$**

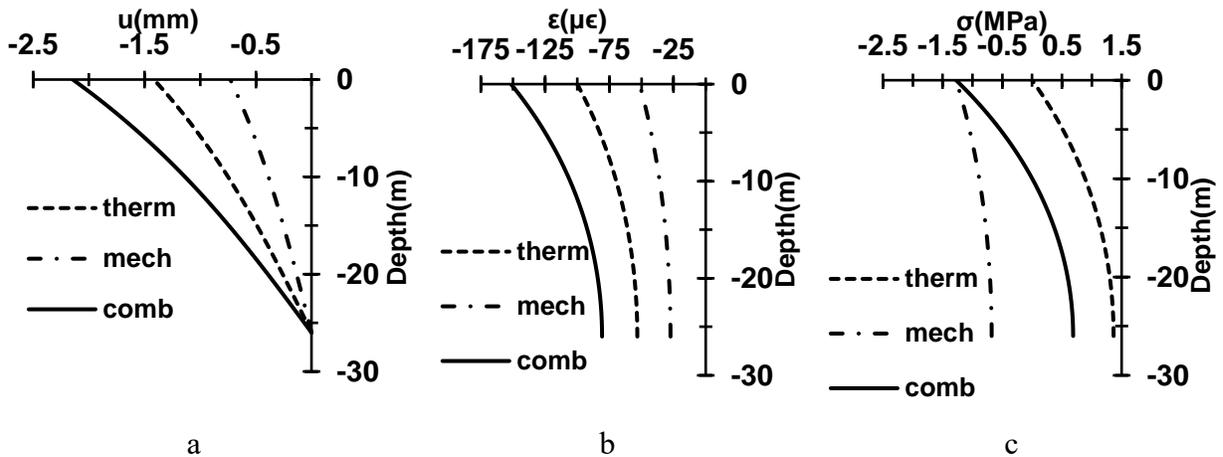

**Figure 3. a) Displacement, (b) strain, and (c) stress in an end bearing energy pile subjected to thermo-mechanical load scenario (*i*), $F = -1000$ kN and $\Delta T = -10$ °C, $k_h = 0$**

Figure 4 shows combined displacement, strain and stress due to load scenario (*ii*) applied to fully floating and end bearing energy piles. It is noted that in this case heating is combined with compressive mechanical load. Thus, in the top portion of the fully floating pile above the thermal null point thermally induced displacement points upward while mechanically induced displacement points downward. This results in the magnitude of combined displacement being smaller than in the case of loading scenario (*i*). Furthermore, in the case of a fully floating pile the combined displacement is nearly zero at the pile head for a given load combination. In the bottom portion of the fully floating pile mechanically and thermally induced displacements are both pointing downward thus resulting in the larger magnitude of the combined displacement at the pile tip than at the pile head. Furthermore, the magnitude of the combined displacement at the pile tip for loading scenario (*ii*) is smaller than the combined displacement at the pile head for the loading scenario (*i*) because of the larger magnitude of mechanically induced displacement in the latter case. Similarly, to displacement response, the magnitude of the combined strain for scenario (*ii*) is smaller than for the scenario (*i*) whereby magnitude of combined strain remains larger for the fully floating pile than for the end bearing pile (Figure 4b). Furthermore, the magnitudes of thermally and mechanically induced stresses in the case of load scenario (*ii*) add up to produce larger magnitude of a combined stress than in case of the load scenario (*i*). As expected, the magnitude of combined stress is larger in the case of end bearing pile than in the case of fully floating pile (Figure 4c).

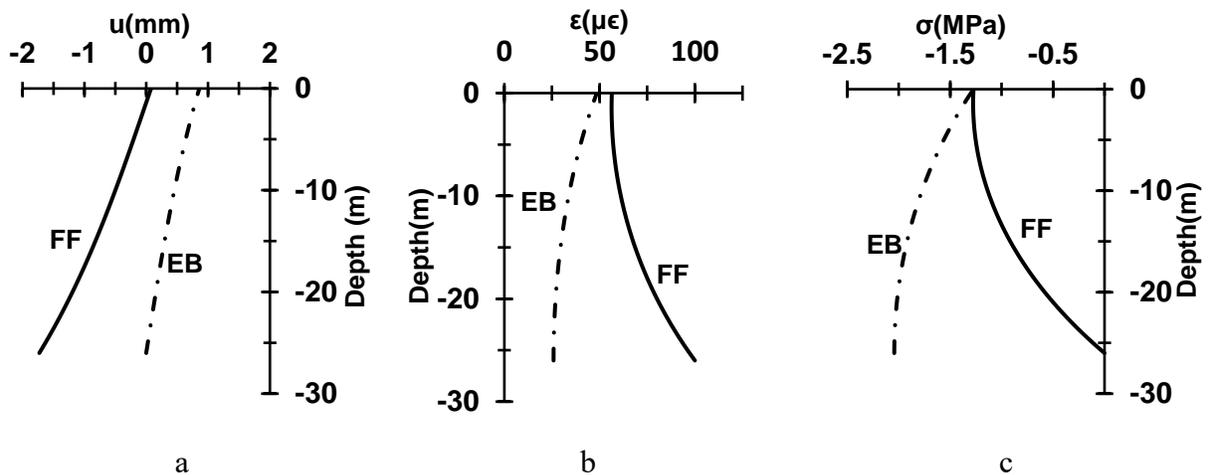

**Figure 4. a) Displacement, (b) strain, and (c) stress in fully floating (FF) and end bearing (EB) energy piles subjected to thermo-mechanical load scenario (*ii*), *F* = -1000 kN and *ΔT*=+10 ᵒC, *$k_h$* =0**

If energy pile that supports a superstructure is thermally loaded the presence of a superstructure will likely impose an additional restraint on the pile displacement. Knellwolf et al. (2011) attached a normal spring having stiffness $k_h$ to the pile head to model this restraint. The effect of $k_h$ is quantified by Eqs. 1 through 3 and 5 through 8 for end bearing and fully floating piles respectively.

In Figure 5 the effect of head restraint on thermo-mechanical displacement of fully floating (FF) and end bearing (EB) energy piles subjected to load scenario (*ii*) is depicted. The selected values of $k_h$ are 0.125 GPa/m and 2 GPa/m as reported by Perić et al. (2017) and Knellwolf et al. (2011) respectively. The thermo-mechanical axial displacement of the end bearing energy pile turns out to be highly restricted by the presence of the head spring. Specifically, its head displacement decreases more than 8 times in case of $k_h = 0.125$ GPa/m as compared Figure 4a and becomes negative for $k_h = 2$ GPa/m. The head of the fully floating energy pile moved even further downward than that of the end bearing pile as the stiffness of the head spring increased.

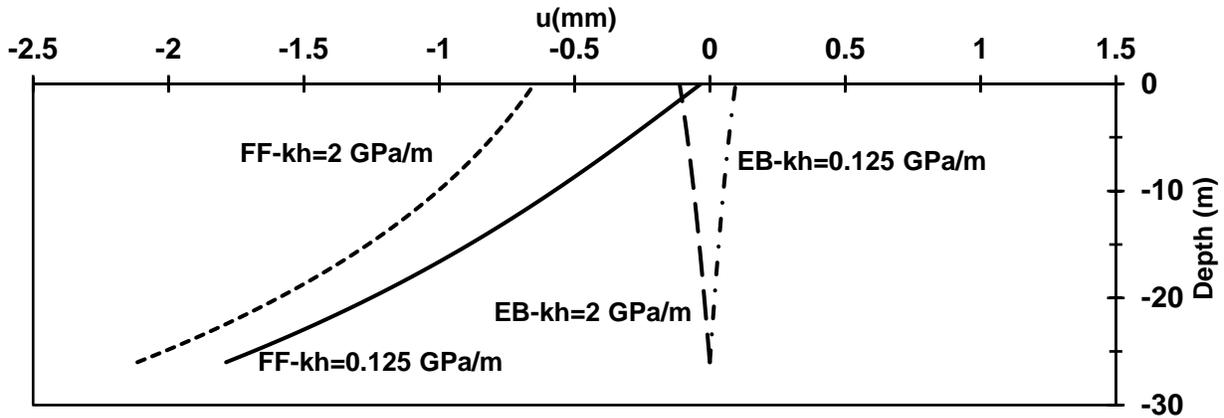

**Figure 5. Thermo-mechanical axial displacement of fully floating (FF) and end bearing (EB) energy piles subjected to load scenario (*ii*), $F = -1000$ kN and $\Delta T = +10$ °C, and $k_h = 0.125$ GPa/m and 2 GPa/m**

Axial strain of energy piles was even more affected by the presence of head spring as compared to axial displacement. Comparing Figure 4b and Figure 6 highlights that the axial strain at the pile head is impacted the most, and head spring shifts the tensile axial strain at the pile head to compressive in both cases of end bearing and fully floating piles. Nevertheless, the absence of restraint at the pile tip provided the freedom for the bottom portion of the fully floating pile to expand. On the contrary, the tip restraint in the end bearing pile forced the strain to remain compressive along the entire length of the pile. Again, the presence of head spring affected the end bearing pile more significantly than the fully floating pile.

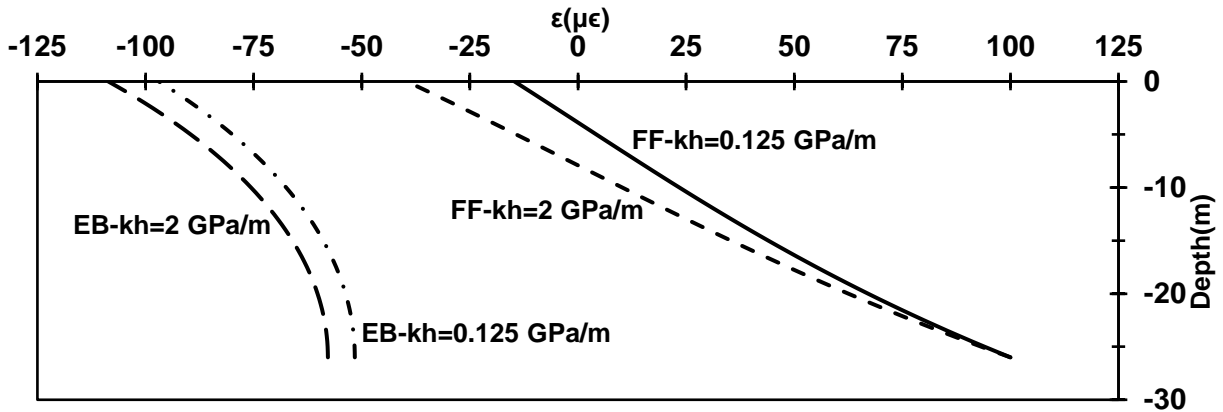

**Figure 6. Thermo-mechanical axial strain of fully floating (FF) and end bearing (EB) energy piles under load scenario (*ii*), $F = -1000$ kN and $\Delta T = +10$ °C, $k_h = 0.125$ GPa/m and $k_h = 2$ GPa/m**

Figure 7 shows that magnitude of compressive thermo-mechanical axial stress increases as the head restraint stiffness ($k_h$) applied at the pile head increases in both, fully floating (FF) and end bearing (EB) energy piles. It is noted that larger values of the thermo-mechanical axial stress were induced in the end bearing pile than in t the fully floating pile. In case of fully floating pile in the magnitude of stress induced by the head restraint decreases with depth and it completely disappears at the pile tip, while the end bearing pile witnesses the compressive stress increase throughout its depth.

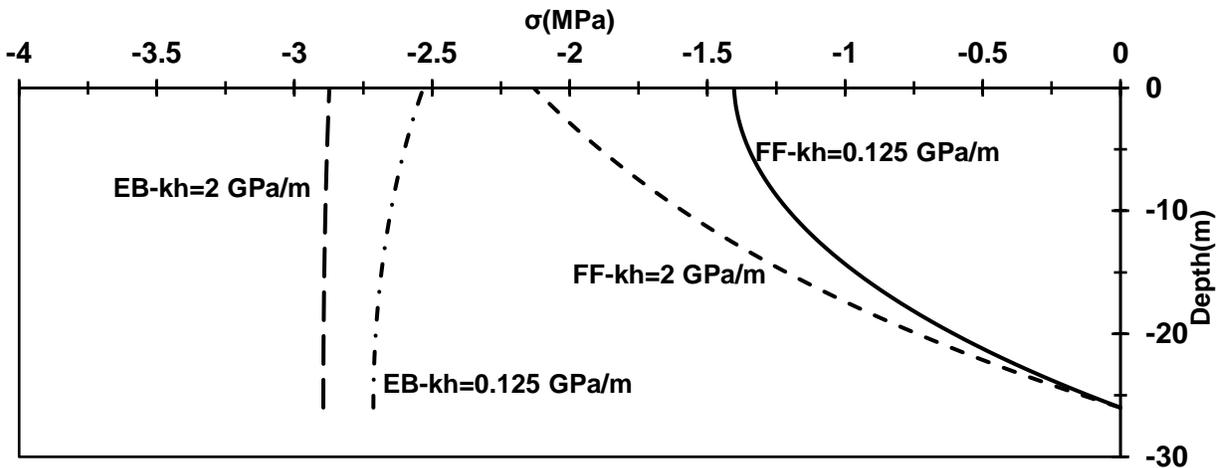

**Figure 7. Thermo-mechanical axial stress of fully floating (FF) and end bearing (EB) energy piles under loading scenario (*ii*), $F$ = -1000 kN and $\Delta T$= +10 ºC, $k_h$= 0.125 GPa/m and $k_h$= 2 GPa/m**

**CONCLUSION**

Effects of end restraints on thermo-mechanical behavior of energy piles were investigated for two load scenarios including combinations of compression and cooling, and compression and heating. End bearing energy piles experience no displacement at the pile tip, which causes a larger magnitude of head displacement in the end bearing than in the fully floating pile. As expected, the end bearing energy pile develops axial stress at its tip while fully floating pile does not. This results in larger development of larger tensile stress in the end bearing pile than in the fully floating pile in case of combined compression and cooling.

Effects of the restraint imposed by interaction of the pile and superstructure were modeled by attaching a normal spring at the pile head. For sufficiently large values of the stiffness of a head spring the thermo-mechanical displacement at the pile head changed the sign, thus changing its direction from upward to downward directed displacement for combined compressive mechanical load and heating. It also increased the compressive stress in the pile. The effect of head restraint on displacement, strain and stress is felt throughout the larger depth of the pile in case of the bearing than in the case of the fully floating pile.